\documentclass[15pt,a4paper]{article}
\usepackage{blindtext}
\usepackage{hyperref}
\usepackage{indentfirst}
\usepackage{misccorr}
\usepackage{graphicx}
\usepackage[all,cmtip]{xy}
\usepackage{amsfonts}
\usepackage{float}
\usepackage{setspace}
\usepackage{eufrak}
\usepackage{wrapfig}
\usepackage[utf8]{inputenc}
\usepackage{csquotes}
\usepackage[T2A]{fontenc}
\usepackage{amsmath,amssymb,amsthm}
\usepackage{cmll}
\usepackage[a4paper,hmargin=2.5cm,vmargin=2.5cm]{geometry}
\usepackage[english]{babel}
\usepackage{cleveref}
\usepackage[toc, page]{appendix}
\RequirePackage[color       = blue!20,
                bordercolor = black,
                textsize    = tiny,
                figwidth    = 0.99\linewidth]{todonotes}

\RequirePackage[backend    = biber,
                sortcites  = true,
                giveninits = true,
                doi        = false,
                isbn       = false,
                url        = false,
                maxnames   = 6,
                style      = numeric]{biblatex}
\DeclareFieldFormat*{title}{#1}
\renewbibmacro{in:}{}

\usetikzlibrary{matrix, arrows, decorations.pathmorphing}
\numberwithin{equation}{section}
\theoremstyle{plain}

\newtheorem*{thmX}{Theorem}

\newtheorem{lemma}[equation]{Lemma}
\newtheorem{theorem}[equation]{Theorem}

\theoremstyle{remark}
\newtheorem{remark}{Remark}
\theoremstyle{definition}

\newcommand{\DM}{\mathbf{DM}}

\newcommand{\SmAff}{\mathrm{SmAff}}

\newcommand{\SH}{\mathbf{SH}}
\newcommand{\HH}{\mathbf{H}}
\newcommand{\Spec}{\operatorname{Spec}}

\newcommand{\pair}{\mathrm{pair}}

\newcommand{\SmAffpair}{\SmAff^\pair}

\newcommand{\X}{X_0}
\newcommand{\Xp}{X_1}
\newcommand{\Z}{Z_0}
\newcommand{\Zp}{Z_1}
\newcommand{\Zcal}{Z^{(1)}_0}
\newcommand{\Zpacal}{Z^{(1)}_1}
\newcommand{\Ztilcal}{\mathcal{Z}}
\newcommand{\Ztil}{Z}
\newcommand{\Zinf}{Z^{(1)}}
\newcommand{\Zinfp}{V}
\newcommand{\Xtil}{X}
\author{
Druzhinin A.E\footnote{
research is supported by PJSC «Gazprom Neft», aggr. ГПН-26/03000/00502/Р}.,
Urazbaev A.A.
}
\title{
On the motivic classification of 
codimension filtration quotients.
}

\begin{document}

\title{Nisnevich equivalences of smooth affine open pairs.}

\maketitle

\section{Introduction}
The Nisnevich topology on the category of algebraic-geometric schemes over a base scheme $B$ takes a place between the Zariski topology and the \'etale topology and allows one to introduce analogues of tubular neighbourhoods of differential submanifolds or topological subspaces in algebraic-geometric geometry, i.e., for closed subschemes $Z$ of a scheme $X$. 
%
Equivalences of the form
\begin{equation}\label{eq:XX-ZsimeqXX-Zprime}X/(X-Z)\simeq X^\prime/(X^\prime-Z^\prime)\end{equation}
allows equivalence
reducing 
the study of motivic invariants of a scheme $X$ to the schemes $X - Z$ and $Z$, together with certain local data: namely, how $Z$ is glued in, i.e., to the normal bundle for smooth $Z$ and $X$
and are useful for the constructions of six functors and duality theorems for stable motivic homotopy categories. 
%
%
%
However, the proofs 
based on deformation to the normal cone or on localization theorems, use not only Nisnevich-local equivalences, but also $\mathbb A^1$-homotopy equivalences.

%
%
We obtain a result for 
affine smooth schemes $X$ and $X^\prime$ and 
arbitrary 
closed subschemes $Z$ and $Z^\prime$ without use of $\mathbb A^1$-homotopy equivalences and only Nisnevich-local equivalences. 

    \begin{thmX}[\Cref{thm:mainthm}]\label{thm:theoremA}
            Let $(\X,\Z), (\Xp,\Zp) \in \SmAffpair_B$ are such that there is the following isomorphism 
        \begin{equation}\label{eq:isoclosub}
            c:\Z\xrightarrow{\cong} \Zp
        \end{equation}
        Then following conditions are equivalent 
        \begin{itemize}
            \item[1.]\label{item:thm1} There is the following isomorphism   
    \begin{equation}\label{eq:isonei}
    c_1:\Zcal\xrightarrow{\cong} \Zpacal
    \end{equation}
     agreed with $\eqref{eq:isoclosub}$.
    \item[2.]\label{item:thm2} There is the following isomorphism
    \begin{equation}\label{eq:isonei2}
        c_2:\Z^{[1]}\xrightarrow{\cong}\Zp^{[1]}
    \end{equation}
     agreed with
     $\eqref{eq:isoclosub}$.
    \item[3.]\label{item:thm3}
    There is the following isomorphism 
    \begin{equation}\label{eq:isonei3}
        c_3:z_0^*\Omega^1_{\X/B}\cong z_1^*\Omega^1_{\Xp/B}
    \end{equation}
    agreed with the isomorphism  $\Omega^1_{\Z/B}\cong\Omega^1_{\Zp/B}$, induced by the isomorphism \eqref{eq:isoclosub}. Here  $z_0:\Z\hookrightarrow\X$, $z_1:\Zp\hookrightarrow\Xp$ denote the embeddings. Hence the following diagram is commutative 
        \begin{equation}\label{eq:isobunZ}
            \xymatrix{
            z_0^*\Omega^1_{\X/B}\ar[d]\ar[r]^{c_3}& 
            z_1^*\Omega^1_{\Xp/B}\ar[d]\\
            \Omega^1_{\Z/B}\ar[r]^{\cong}&\Omega^1_{\Zp/B}
            }
        \end{equation}
    \item[4.]\label{item:thm4}
    There exists a scheme $\Xtil$, a closed subscheme $Z$ and etale morphisms $p_0$, $p_1$ such that the following diagram is commutative
    \begin{equation}\label{eq:XZzXZXZo}
        \xymatrix{
          & \Xtil\ar[rd]^{p_1} \ar[ld]_{p_0} & \\
          \X & Z\ar[rd]^{\cong}\ar[ld]_{\cong}\ar@{^(->}[u] & \Xp\\
          \Z\ar@{^(->}[u]& &\Zp\ar@{^(->}[u]
        }
    \end{equation}
    \end{itemize}
    \end{thmX}

The diagram \eqref{eq:XZzXZXZo} 
provides a pair of Nisnevich squares 
between 
$X_0/(X_0-Z_0)$ and $X_1/(X_1-Z_1)$
that define 
any of the above points implies
a Nisnevich local equivalences
$X_0/(X_0-Z_0)\simeq X_1/(X_1-Z_1)$.

\begin{remark}
The result for local schemes was proven in our earlier preprint \cite{DruUr-LocalEssSmPairsNisEq} by 
a shorter and more elementary argument. The present proof is not based on the result for local schemes,
and we have doubts about such an argument that would deduce the result for given schemes from the results for the local ones. 
The proof is obtained by a direct global construction of the diagram \eqref{eq:XZzXZXZo}, which uses affineness of the schemes.
\end{remark}


    \subsection{Agreements and notation}
\begin{enumerate}
\item
$X-Z$ and $X\setminus U$ are open and reduced closed complements, respectively, to the closed subscheme $Z$ and the open subscheme $U$ in $X$.
\item \label{notation:categoriesSH(B)DM(B)}
$\mathbf{SH}(B)$ ---
the stable motivic homotopy category of Voevodsky over $B$, in the sense of
\cite{MV,Voe98,Jardine-spt},
\cite[Appendix C]{Hoyois_SHbaseschemeBchN}. 
$\mathbf{DM}(B)$ ---
the category of motives of Voevodsky over $B$, in the sense of
\cite{Voe-motives,Cisinski-Deglise-Triangmixedmotives}.
    \item 
    $\SmAffpair_B$ -- the category of pairs $(X,Z)$, where $X$ is a smooth scheme over a Noetherian separated scheme $B$, and $Z$ is a closed not necessarily smooth subscheme of $X$.
    \item 
    $Z^{(1)}_X$ --  the first infinitesimal neighbourhood of the closed subscheme $Z$ in $X$, defined by $Z^{(1)}_X = Z(I^2_X(Z))$. If the scheme $X$ is understood, we will use the abbreviated notation $Z^{(1)}$.
    \item 
    $\Omega^1_{X/B}$ -- the module of Kähler differentials of $X$ over $B$.
    \item
    $d_ZX:I/I^2\to z^*\Omega^1_{X/B}$ -- the morphism, where $Z$ is a closed subscheme in $X$ and $z:Z\hookrightarrow X$ denotes a closed embedding(see \cref{lm:etalcri}).
    \item 
    $Z^{[1]}_X$ -- the closed subscheme of $X$ that contains the closed subscheme $Z$ and is equal to  $\Spec(A/J)$, where $X=\Spec(A)$, $Z=\Spec(A/I)$ and $J$ is inverse image of $ker(d_ZX)$ under $I\twoheadrightarrow I/I^2$.
    \end{enumerate}

    \section{Affine pairs}
        Let  $(\X,\Z), (\Xp,\Zp) \in \SmAffpair_B$ be such that there is the following isomorphism  \begin{equation}\label{eq:isodef}
    c:\Z\xrightarrow{\cong} \Zp
    \end{equation} 
    Let $\Ztil$ denote the graph of the isomorphism \eqref{eq:isodef}, and let $\Ztilcal$ be a closed subscheme of $\X\times_B\Xp$ containing $\Ztil$ such that the projections $p_0$, $p_1$ induce isomorphisms from $\Ztilcal$ onto its images $\Ztilcal_0$ and $\Ztilcal_1$, respectively.
  Let   
    \begin{equation}\label{eq:defV}
    V:=\Ztilcal^{(1)}_{\Ztilcal_0\times_B\Xp}.
    \end{equation}
    Note that $\Ztilcal$, $\Ztilcal_0\times_B\Xp$ are smooth over $\Ztilcal_0$. Consider the scheme  $N$ corresponding to the normal bundle  $N_{\Ztilcal \hookrightarrow \Ztilcal_0 \times_B\Xp}$ and the commutative diagram 
    \begin{equation}\label{}
        \xymatrix{
        N\ar[d]&& N^{\times 2}_{\Ztilcal}\ar[ll]\ar@<-1ex>[d]\\
        \Ztilcal\ar@<-1ex>[rr]_{0}&& N\ar@<-1ex>[ll]\ar@<-1ex>[u]_{\Delta}
        }
        \end{equation}
        where $\Delta$  is the diagonal embedding and $0$ is the zero section.\par
        Since $\Ztilcal$, $\Ztilcal_0\times_B\Xp$ are smooth over $\Ztilcal_0$, there is also a morphism $\Ztilcal^{(1)}\xrightarrow{p_0}\Ztilcal_0\cong \Ztilcal$. Then by  $\cref{prop:norbun}$, there is the isomorphism  \begin{equation}\label{eq:is_norbun}
            \Zinfp \cong Z(I^2_{N_{\Ztilcal \hookrightarrow \Ztilcal_0 \times_B\Xp}}(0)),
        \end{equation} hence there is a closed embedding  $\Zinfp\hookrightarrow N$. Denote by  $\Delta^\prime:=\Delta\big|_{\Zinfp}$.
        \begin{remark}
        Note that there is a closed embedding $i:\Ztilcal\hookrightarrow\Zinfp$; also $\Delta^\prime\circ i=0$, because the isomorphism \eqref{eq:isoinfnei} is agreed with the isomorphism $\Ztilcal\cong 0$. 
        \end{remark}
        Consider the Cartesian square
    \begin{equation}
        \xymatrix{
        N^\prime\ar[d]\ar[rr]& & N^{\times2}_{\Ztilcal}\ar[d]\\
        V\ar[rr]& & N
        }
    \end{equation}
    and consider the Cartesian diagram
     \begin{equation}\label{eq:DiagN}
        \xymatrix{
        N^\prime\ar[rr]& & N^{\times2}_{\Ztilcal}\\
        V\ar[rr]\ar[u]^{\Delta^\prime}& & N\ar[u]^{\Delta}
        }
    \end{equation}
    where $\Delta^{\prime}:=\Delta\big|_{V}$.
    \begin{lemma}\label{lm:exbun}      
    Note that there is a closed embedding $V\hookrightarrow X_0\times_BX_1$, where  $V$, as in  $\eqref{eq:defV}$. Let  $U:=(\X\times_B\Xp)^h_{V}$, then there exists a vector bundle $M$ over $U$ together with a section $s:U\to M$ such that  $M\times_U\Zinfp\cong N^\prime$, where $N^\prime$, as in $\eqref{eq:DiagN}$ and the following diagram is commutative 
        \begin{equation}\label{diaf:conbun}
            \xymatrix{
            N^\prime\ar[rr]&& M\\
            V\ar[rr]\ar[u]^{\Delta^\prime} && U\ar[u]_s
            }
        \end{equation}
    \end{lemma}
    \begin{proof}
        Since our schemes are affine, we may assume that  $U=\Spec(R)$ and $V=\Spec(R/I)$. Then, by \cref{lm:0D4A}, there exists a vector bundle $M$. And by  \cref{lm:0H74} applied to   $X:=U$, $Z:=V$, $Y:=M$, $B:=U$, there exists an extension of the morphism $V\to M$, and since $s$ is a morphism of schemes over $U$, it is a section.  
    \end{proof}  
    \begin{theorem}\label[theorem]{thm:mainthm}
        Let $(\X,\Z), (\Xp,\Zp) \in \SmAffpair_B$ are such that there is the following isomorphism 
        \begin{equation}\label{eq:isoclosub}
            c:\Z\xrightarrow{\cong} \Zp
        \end{equation}
        Then following conditions are equivalent 
        \begin{itemize}
            \item[1.]\label{item:thm1} There is the following isomorphism   
    \begin{equation}\label{eq:isonei}
    c_1:\Zcal\xrightarrow{\cong} \Zpacal
    \end{equation}
     agreed with $\eqref{eq:isoclosub}$.
    \item[2.]\label{item:thm2} There is the following isomorphism
    \begin{equation}\label{eq:isonei2}
        c_2:\Z^{[1]}\xrightarrow{\cong}\Zp^{[1]}
    \end{equation}
     agreed with
     $\eqref{eq:isoclosub}$.
    \item[3.]\label{item:thm3}
    There is the following isomorphism 
    \begin{equation}\label{eq:isonei3}
        c_3:z_0^*\Omega^1_{\X/B}\cong z_1^*\Omega^1_{\Xp/B}
    \end{equation}
    agreed with the isomorphism  $\Omega^1_{\Z/B}\cong\Omega^1_{\Zp/B}$, induced by the isomorphism \eqref{eq:isoclosub}. Here  $z_0:\Z\hookrightarrow\X$, $z_1:\Zp\hookrightarrow\Xp$ denote the embeddings. Hence the following diagram is commutative 
        \begin{equation}\label{eq:isobunZ}
            \xymatrix{
            z_0^*\Omega^1_{\X/B}\ar[d]\ar[r]^{c_3}& 
            z_1^*\Omega^1_{\Xp/B}\ar[d]\\
            \Omega^1_{\Z/B}\ar[r]^{\cong}&\Omega^1_{\Zp/B}
            }
        \end{equation}
    \item[4.]\label{item:thm4}
    There exists a scheme $\Xtil$, a closed subscheme $Z$ and etale morphisms $p_0$, $p_1$ such that the following diagram is commutative
    \begin{equation}
        \xymatrix{
          & \Xtil\ar[rd]^{p_1} \ar[ld]_{p_0} & \\
          \X & Z\ar[rd]^{\cong}\ar[ld]_{\cong}\ar@{^(->}[u] & \Xp\\
          \Z\ar@{^(->}[u]& &\Zp\ar@{^(->}[u]
        }
    \end{equation}
    \end{itemize}
    \begin{proof} 
    \par\noindent\textbf{ $(1)\Rightarrow (4)$}\par
    Let us define $\Xtil:=Z(s)$, where $s$ is as in \eqref{diaf:conbun} for $V$ as in \eqref{eq:defV} with  $\Ztilcal$ being the graph of the isomorphism $\eqref{eq:isonei}$ and $\Ztil$ being the graph of the isomorphism  $\eqref{eq:isoclosub}$. It is obvious that  $\Ztilcal$ maps isomorphically onto its image under the projection. By the construction $\Ztil\hookrightarrow\Xtil$. Morphisms $p_0$, $p_1$ are a compositions of the embedding and the following projection. Isomorphism onto $\Ztil$ is obvious. It remains only to check that the morphisms are etale. \\
    Recall that $V=\Ztilcal^{(1)}_{\Zcal\times_B\Xp}$ also  \[V\times_{\Zcal\times_B\Xp}\Xtil\cong \Ztilcal,\] since, by the commutativity of the diagram\eqref{diaf:conbun} follows that $V\times_{\Zcal\times_B\Xp}\Xtil\cong Z(\Delta^\prime)$. Hence \begin{equation}\label{eq:lemiso}
             \Ztilcal^{(1)}_{\X\times_B\Xp}\times_{\X\times_B\Xp}\Xtil\times_{\X}\Zcal\cong \Ztilcal
         \end{equation}  
         By the  \cref{lm:heneq}, we may assume that $\X\cong(\X)^{h}_{\Z}$, $\Xp\cong(\Xp)^{h}_{\Zp}$, also $(\X\times_B\Xp)_\Ztil^h\cong(\X\times_B\Xp)_V^h$. Therefore, we can apply  \cref{lm:eqprop} for $\Ztilcal$, obtaining that the morphism $p_0$ is etale. We also note that $\Zinf\subseteq\Zpacal\times_B\Zcal$; hence  \[\Zp\times_{\Xp}\Zinf\subseteq\Zp\times_{\Xp}\Zinf\times_{\X}\Zcal\] 
          and
          \[\Zp\times_{X_1}\Zinf\times_{\X\times_B\Xp}\Xtil\subseteq\Zp\times_{X_1}\Zinf\times_{\X\times_B\Xp}\Xtil\times_{\X}\Zcal\cong\Zp\times_{\Xp}\Ztil\cong\Ztil
          \]
          The first isomorphism holds by etaleness of $p_0$.
          The opposite inclusion $ \Ztil\subseteq \Zp\times_{X_1}\Zinf\times_{\X\times_B\Xp}\Xtil$is obvious. Hence, by \cref{lm:eqprop} the morphism $p_1$ is also etale.
          \par\noindent\textbf{ $(4)\Rightarrow (3)$}\par

       From the etaleness of the morphism $p_0:X\to X_0$, it follows that  $z^*\Omega^1_{X/B}\cong z_0^*\Omega_{X_0/B}$. Similarly,  $z^*\Omega^1_{X/B}\cong z_1^*\Omega_{X_1/B}$.
          \par\noindent\textbf{ $(3)\Rightarrow (2)$}\par
        Let us consider the subbundle $T\hookrightarrow z^*\Omega^1_{\X\times_B\Xp/B}\cong z_0^*\Omega^1_{\X/B}\oplus z_1^*\Omega^1_{\Xp/B}$ equals to the graph of the isomorphism  $c_3$ from \eqref{eq:isobunZ}. Let us consider the ideal $J:=\{f\in\mathcal{O}_{\X\times_B\Xp}\big|f|_{\Ztil}=0, d_{\Ztil}f\in T\}$. Let us define $\Ztilcal:=Z(J)$, let us check that it maps isomorphically onto $\Tilde{\Z}$, $\Tilde{\Zp}$ under the projections. We consider the following morphism \begin{equation}
        \sigma:z^*\Omega^1_{\Xtil/B}\cong z_0^*\Omega^1_{\X/B}\oplus z_1^*\Omega^1_{\Xp/B}\to z_1^*\Omega^1_{\X/B}
        \end{equation} defined as follows $a\oplus b\mapsto b-c_3(a)$. Notice that $\sigma^\prime$ vanishes on $J$. Consider the following commutative diagram
        \begin{equation}
            \xymatrix{
            I\ar[d]\ar[rd]_{\sigma^\prime}\ar[r]^-{d_Z}&z^*\Omega^1_{\X\times_B\Xp/B}\ar[d]^{\sigma}\ar[r]& 
            \Omega^1_{\Ztil/B}\ar[r]&0\\
            I/J\ar[r]&z_0^*\Omega^1_{\X/B}\ar[r]^{i_0}&\Omega^1_{\Z/B}\ar[r]&0
            }
        \end{equation}
        here $I$, $I_0$ and  $I_1$ are ideals corresponding to $\Ztil$, $\Z$ and $\Zp$. The upper sequence is the one from \eqref{eq:exseqdif}. We shall prove that the lower sequence is short exact. 
        Exactness at $\Omega^1_{\Z/B}$  follows from the exactness of \eqref{eq:exseqdif} for $\Z\hookrightarrow \X$.Exactness at $\Z^*\Omega^1_{\X/B}$ by the commutativity of the following diagram
        \begin{equation}
            \xymatrix{
            I_1\ar[d]_{p_1^*}\ar[r]&z_1^*\Omega^1_{\X/B}\ar[d]^{d_1}\\
            I\ar[r]&z^*\Omega^1_{\X\times_B\Xp/B}\cong z_0^*\Omega^1_{\X/B}\oplus z_1^*\Omega^1_{\Xp/B}\\
            I_0\ar[u]^{p_0^*}\ar[r]&z_0^*\Omega^1_{\X/B}\ar[u]_{d_0}
            }
        \end{equation}
        because $i_0\circ\sigma(a\oplus b)=i_0\circ d_0(a)-i_0\circ d_0(c_3(b))$, hence $i_0\circ\sigma^\prime=0$. Also $\sigma\circ d_0=id$.
        Since $J$ coincides with $ker(\sigma^\prime)$, we have that $I/J\cong I_0/J_0$, where $J_0$ is a kernel of the map $I_0\to\Z^*\Omega^1_{\X/B}$. Also $\Z\cong\Zp$, therefore $\Z^{[1]}\cong\Zp^{[1]}$.
        \par\noindent\textbf{ $(2)\Rightarrow (4)$}
        \par
        Consider $\Xtil:=Z(s)$, where $s$ is as in \eqref{diaf:conbun} for $V$ as in \eqref{eq:defV}, with $\Ztilcal$ being the graph of the isomorphism \eqref{eq:isonei2} and $\Ztil$ being the graph of the isomorphism \eqref{eq:isoclosub}.
        $\eqref{eq:isoclosub}$.
         The morphism $p_0:X\to\X$ induces the following commutative diagram 
        \begin{equation}
            \xymatrix{
            I/J\ar[r]&z^*\Omega^1_{X/B}\ar[r]&\Omega^1_{Z/B}\ar[r]&0\\
            I_0/J_0\ar[u]^{\cong}\ar[r]&z_0^*\Omega^1_{\X/B}\ar[u]\ar[r]&\Omega^1_{\Z/B}\ar[u]^{\cong}\ar[r]&0
            }
        \end{equation}
         By the five-lemma, $z_0^*\Omega^1_{\X/B}\cong z^*\Omega^1_{X/B}$.Then by \cref{lm:etalcri} the morphism $p_0$ is etale.The proof for $p_1$ is analogous.
        \par\noindent\textbf{ $(4)\Rightarrow (1)$}\par
        By \cref{lm:heneq}, we have the isomorphism $(\X)^h_{\Z}\cong(\Xp)^h_{\Zp}$, This implies the isomorphism \eqref{eq:isonei}.
        
    \end{proof}
    \end{theorem}
     \begin{remark}\label{cor:mainthm}
    Let $(\X,\Z), (\Xp,\Zp) \in \SmAffpair_B$ be such that there is the following isomorphism
    \begin{equation}\label{eq:isoclosub111}
            c:\Z\xrightarrow{\cong} \Zp
        \end{equation}
        and one of the conditions (1)--(4) of \Cref{thm:mainthm} is satisfied. Then there is the isomorphism of Nisnevich presheaves
        \begin{equation}
            \X/(\X-\Z)\cong\Xp/(\Xp-\Zp)
        \end{equation}
        and hence an isomorphism in  categories $\HH(B)$, $\SH(B)$ and $\DM(B)$
    \end{remark}
    \section{Appendix A}            
    \begin{lemma}\label{lm:eqprop}
          Let $(\X,\Z), (\Xp,\Zp) \in \SmAffpair_B$ be such that $\X\cong(\X)^{h}_{\Z}$, $\Xp\cong(\Xp)^{h}_{\Zp}$ and suppose there is an isomorphism $\Z\xrightarrow{\cong} \Zp$. Let $\Ztil$ be the graph of this isomorphism. Consider $\Xtil$ a closed subscheme in $U:=(\X\times_B\Xp)_\Ztil^h$ containing $\Ztil$. Let $p_0:\Xtil\to \X$ be the composition of the morphism $\Xtil\hookrightarrow U\rightarrow\X\times_B\Xp$  and the projection $\X\times_B\Xp\to \X$. We identify $\Ztil$ with its image under the morphism   $f:\Ztil\to\Xtil\times_{X_0}Z_0$. 
            Then the following statements are equivalent
            \begin{itemize}
               
                \item[1.] 
                $Z$  is both open and closed in $\Xtil\times_{X_0}Z_0$.
                 \item[2.] Morphism $p_0^\prime:\Xtil\times_{X_0}Z_0\to Z_0$ is etale over $Z$.
                \item[3.]
                Morphism $p_0:\Xtil\to X_0$ is etale over $Z$.
                \item[4.]
                Morphism  $p_0:\Xtil\to X_0$ is etale.
                \item[5.]
                 There is the isomorphism $\Ztil\xrightarrow{\cong}\Zinf\times_{Z_0\times_B\Xp}(\Xtil\times_{X_0}\Z)$. 
                \end{itemize}
                \begin{proof}
                \par\noindent\textbf{ $(1)\Rightarrow (2)$}\par
                     In this case, the morphism $\Ztil\hookrightarrow{} \Xtil\times_{\X}\Z$ is an open embedding, hence an etale morphism. Therefore, in the following commutative diagram, 2 of 3 morphisms are etale.
                     \begin{equation}
                         \xymatrix{
                         Z\ar[r]^-{et}\ar[d]^{\cong}& \Xtil\times_{X_0}Z_0\ar[ld]^{p_0^\prime}\\
                         \Z
                         }
                     \end{equation}
                     and morphism  $p_0^{\prime}$ are etale, because etale morphisms satisfy the 2 of 3 property. 
                     \\
            \par\noindent\textbf{ $(2)\Rightarrow (3)$}\par
                    Since $p_0:X\to X_0$ is locally of finite presentation  by (2) from  \cite[\href{https://stacks.math.columbia.edu/tag/02GU}{Lemma 02GU}]{stacks-project}
                    it suffices to show that the morphism of local rings $\mathcal{O}_{X,x}\to\mathcal{O}_{X_0,x_0}$ is flat, where $x\in Z$,  $x_0:=p_0(x)$. Let us consider the following commutative diagram of rings
                    \begin{equation}
                        \xymatrix{
                        B\ar[r]&B\otimes_{A_0}A_0/J_0\\
                        A_0\ar[r]\ar[u]&A_0/J_0
                        \ar[u]}
                    \end{equation}
                    where $X=\Spec(B), X_0=\Spec(A_0), x_0=\Spec(A_0/J_0)$
                    such that $A_0/J_0\to B\otimes_{A_0}A/J_0$ is a flat morphism of rings. Consider an exact sequence of $A_0$-modules
                    \begin{equation}
                        \xymatrix{
                        0\ar[r]&M^\prime\ar[r]&M\ar[r]&M^{\prime\prime}\ar[r]&0.
                        }
                    \end{equation}
                    Recall that $A_0\to A_0/J_0$ defines a faithfully flat module, since $\X\cong(\X)^{h}_{\Z}$, and that the composition of flat morphisms is flat. Hence the following sequence is exact. 
                    \begin{equation}
                        \xymatrix{
                            0\ar[r]&M^\prime\otimes_{A_0} B\otimes_{A_0}A_0/J_0\ar[r]&M\otimes_{A_0} B\otimes_{A_0}A_0/J_0\ar[r]&M^{\prime\prime}\otimes_{A_0} B\otimes_{A_0}A_0/J_0\ar[r]&0.
                        }
                    \end{equation} 
                    and since the implication in the opposite direction holds for $A_0\to A_0/J_0$, the sequence is exact
                    \begin{equation}
                        \xymatrix{
                        0\ar[r]&M^\prime\otimes_{A_0} B\ar[r]&M\otimes_{A_0} B\ar[r]&M^{\prime\prime}\otimes_{A_0} B\ar[r]&0.
                        }
                    \end{equation}
                \par\noindent\textbf{ $(3)\Rightarrow (4)$}\par
                    The morphism $p_0$ is etale in some open neighbourhood of $\Ztil$ and since $X$ is a closed subscheme of $(\X\times_B\Xp)^h_{\Ztil}$,  this neighbourhood coincides with $\Xtil$.\\
                    \par\noindent\textbf{ $(4)\Rightarrow (5)$}\par
                     It is enough to prove that \[\Xtil\times_{\X}\Ztil\cong\Ztil\sqcup W,\] hence  \[\Zinf\times_{Z_0\times_B\Xp}\Xtil\times_{X_0}\Z\cong \Zinf\times_{\Z\times_B\Xp}\Ztil\sqcup \Zinf\times_{\Z\times_B\Xp}W\cong \Ztil\sqcup \emptyset.\] It is obvious that  $f(\Ztilcal)$ is closed.
                    Since etaleness is stable under base change and is an open property, the morphism $p^\prime_0$ is etale on $W$ an open neighbourhood of $Z$ in $\Xtil\times_{X_0}Z_0$.  The morphism $p_0^\prime$ is unframed, hence the diagonal morphism \[\Delta:W\to W\times_{Z_0}W\] is an open embedding. Morphism $p_0^\prime$ has the section \[i:Z_0\xrightarrow{\cong}Z\xrightarrow{f}\Xtil\times_{X_0}Z_0.\] Let us consider the morphism 
                    \[(i\circ p_0^\prime,id):W\to W\times_{\Zcal}W.\] 
                    Hence the diagonal is open and $(i\circ p_0^\prime,id)^{-1}$ is continuous; then \[(i\circ p_0^\prime,id)^{-1}(\Delta_{W\times_{\Zcal}W})=f(\Ztilcal)\] is open.\\
                \par\noindent\textbf{ $(5)\Rightarrow (1)$}\par
                    Let us restrict to the scheme  $\Z\times_B\Xp$. Let $I$ be the ideal defining $\Ztil$ and let $J$ be the ideal defining $\Xtil\times_{\X}Z_0$. Then from $\Ztil \xrightarrow{\cong}V\times_{Z_0\times_B\Xp}\Xtil\times_{X_0}\Z$ it follows that $J+I^2=I$ and \begin{equation}
                             I(I/J)=I/J.
                         \end{equation} 
                         Let $\eta$ define a common point of the scheme $\Ztil$ in $\Z\times_B\Xp$. The ideal $I_\eta$ lies in the Jacobson radical of the local ring, since $\eta\in\Ztil$.  By Nakayama's lemma $(I/J)_\eta=0$, which means that \[\Xtil\times_{\X}\Z\times_{\Z\times_B\Xp}U=\Ztil.\] Hence $\Ztil$ is both open and closed in $\Xtil\times_{\X}\Z$.
                \end{proof}
            \end{lemma}
    \begin{lemma}\label{lm:etalcri}
        Let $(X,Z)\in\SmAffpair$, $X=\Spec(A)$, $Z=\Spec(A/I)$. Let denote by $z:Z\hookrightarrow X$ the embedding. Then there is the following exact sequence         \begin{equation}\label{eq:exseqdif}
            \xymatrix
            {
                I/I^2\ar[r]^-{d_zX}&z^*\Omega^1_{X/B}\ar[r]&\Omega^1_{Z/B}\ar[r]&0
            }
        \end{equation}
        where $z^*\Omega^1_{X/B}=\Omega^1_{X/B}\otimes_A A/I$.
        \begin{proof}
            Since the schemes are affine, the lemma is equivalent to \cite[Chapter VI,Theorem (1.8)]{altman2006introduction}
        \end{proof}
    \end{lemma}
    \begin{lemma}
        Let $(X,Z), (\X, \Z)\in\SmAffpair$ and the morphism $f:X\to \X$ defines the isomorphism $Z\cong\Z$ and the isomorphism $z_0^*\Omega^1_{\X/B}\cong z^*\Omega^1_{X/B}$,  agreed with $\Omega^1_{Z_0/B}\cong \Omega^1_{Z/B}$. Also $X_0=(X_0)^h_{Z_0}$. Then $f$ is etale.
        \begin{proof}
            By (4) from \cite[\href{https://stacks.math.columbia.edu/tag/02GU}{Lemma 02GU}]{stacks-project} it is enough to prove flatness of the morphism of local rings for every  $x\in X$  $\mathcal{O}_{X,x}\to\mathcal{O}_{X_0,x_0}$ and that $\Omega_{X/X_0,x}=0$ where  $x_0=f(x)$. The first follows from the fact that  $X_0=(X_0)^h_{Z_0}$, and the second from the isomorphism $z_0^*\Omega^1_{\X/B}\cong z^*\Omega^1_{X/B}$.  
        \end{proof}
    \end{lemma}
        \begin{lemma}\label{prop:norbun}
        Let $(X,Z)\in \SmAffpair_B$ with $Z$ closed in $X$ such that $X=\Spec(A)$, $Z=\Spec(A/I)$ and it is a quasi-regular immersion. Then 
        \begin{equation}\label{eq:isoinfnei}
           Z^{(1)}_{X}\cong 0^{(1)}_{N_ZX},
        \end{equation} 
        where $0$ is the zero locus in the normal bundle $N_{Z}X$. Moreover, this isomorphism is compatible with the isomorphism $Z\cong 0$.
        \begin{proof}
        According to  \href{https://stacks.math.columbia.edu/tag/063M}{Lemma 063M} for a quasi-regular embedding, the following equality holds $N_ZX = \underline{\mathop{\mathrm{Spec}}}_ Z(\bigoplus _{n \geq 0} I^ n/I^{n + 1})$. Then \[0^{(1)}_{N_ZX}\cong\Spec(A/I\bigoplus I/I^2),\] also there is the following short exact sequence of modules $A/I$ 
        \begin{equation}
            \xymatrix{
            0\ar[r]&I/I^2\ar[r]&A/I^2\ar[r]&A/I\ar[r]&0
            }
        \end{equation}
         Therefore $A/I^2\cong A/I\bigoplus I/I^2$ as modules.  Moreover, the multiplicative structures on these modules coincide, hence $\Spec(A/I^2)\cong\Spec((A/I)\oplus I/I^2)$.
        \end{proof}
    \end{lemma}
    \begin{lemma}\label{lm:0D4A}[\href{https://stacks.math.columbia.edu/tag/0D4A}{Lemma 0D4A}]
        Let (X,Z) is a henselian pair of affine schemes.  Let (X,Z) is henselian pair of affine schemes. The map $E \to E \times_X Z$ induces a bijection between the sets of isomorphism classes of finite-dimensional vector bundles over $X$ and finite-dimensional vector bundles over $Z$. In particular, every finite-dimensional vector bundle $E^\prime$ over $Z$ is isomorphic to $E \times_X Z$ for some finite-dimensional vector bundle $E$ over $X$
        \begin{equation}
            \xymatrix{
            E^\prime\ar[r]\ar[d]& E\ar[d]\\
            Z\ar[r]&X
            }
        \end{equation}
        \end{lemma}
        \begin{lemma}\label{lm:0H74}[\href{https://stacks.math.columbia.edu/tag/0H74}{Lemma 0H74}]
        Let $Y$ be a smooth affine scheme over a base $B$, and let $X$ be an affine scheme over $B$ such that $(X,Z)$ is a Henselian pair. Then any morphism $Z \to Y$ can be lifted to a morphism $X \to Y$.
        \begin{equation}
            \xymatrix{
            Y\\
            Z\ar[u]\ar[r]&X\ar@{-->}[ul]
            }
        \end{equation}
    \end{lemma}
    \begin{lemma}\label{lm:heneq}
        Let $(\X,\Z), (\Xp,\Zp)\in \SmAffpair_B$ such that $\Z\cong\Zp$. Then following conditions are equivalent
        \begin{itemize}
            \item[1.] There exists $(X,Z)\in\SmAffpair_B$  together with etale morphisms $p_0:X\to\X$, $p_1:X\to\Xp$ such that $p_0\big|_{Z}$, $p_1\big|_{Z}$ are isomorphisms.
            \item[2.] $(\X)^h_{\Z}\cong(\Xp)^h_{\Zp}$
        \end{itemize}
        \begin{proof}
            \par\noindent\textbf{ $(1)\Rightarrow (2)$}\par
           Recall that $(X_0)^h_{Z_0} = \varprojlim U_i$, where each $U_i$ is a scheme with a closed subscheme $Z_i$ such that there exists an étale morphism $p_i: U_i \to X_0$ inducing an isomorphism $Z_i \xrightarrow{\cong} Z_0$. Note that $X$ lies in this filtered category. Hence $(\X)^h_{\Z}\cong(X)^h_Z$. Similarly $(\Xp)^h_{\Zp}\cong(X)^h_Z$.
            \par\noindent\textbf{ $(2)\Rightarrow (1)$}\par
            Since $X_0$, $X_1$ are of finite presentation, then by the lemma \cite[\href{https://stacks.math.columbia.edu/tag/00QO}{Tag 00QO}]{stacks-project} there exist indexes $i$ and  $j$ such that $U_i\cong V_j$, where $U_i$ is neighbourhood of $Z_0$, and $V_j$ is neighbourhood of  $Z_1$. Then let $X := U_i$. The etale morphism to $X_0$ is, by definition, the given morphism $U_i \to X_0$, while the morphism to $X_1$ is induced by the isomorphism with $V_j$. 
        \end{proof}
    \end{lemma}
    
\printbibliography

\par\noindent
Andrei Druzhinin; \\
Chebyshev Laboratory, St. Petersburg State University, 14th Line V.O., 29, Saint Petersburg 199178 Russia;\\ 
St. Petersburg Department of Steklov Mathematical Institute 
of Russian Academy of Sciences, Russia;\\
a.druzhinin@spbu.ru,
andrei.druzh@gmail.com.

\par\noindent
Askar Urazbaev; \\
St. Petersburg Department of Steklov Mathematical Institute 
of Russian Academy of Sciences, Russia;\\
St. Petersburg State University, Department of Mathematics and Computer Sciences, 14th Line V.O., 29, Saint Petersburg 199178 Russia.\\
st087096@student.spbu.ru, Urazaskar@yandex.ru

\end{document}